\documentclass[a4paper,11pt,pointlessnumbers]{scrartcl}
\usepackage{amsmath}
\usepackage{amsfonts}
\usepackage{amssymb,mathrsfs,wasysym}
\usepackage{t1enc}
\usepackage[latin1]{inputenc}
\usepackage[greek.ancient,english]{babel}
\usepackage{natbib}
\usepackage{amssymb}
\usepackage[final]{pdfpages}
\usepackage{caption}
\usepackage{url}

\usepackage[colorlinks=false, linktocpage=true,backref,pagebackref,hidelinks]{hyperref}
\begin{document}

\title{\bf \large Mathematical modernity, goal or problem? The opposing views of Felix Hausdorff and Hermann Weyl}
\author{\normalsize Erhard Scholz\footnote{University of Wuppertal, Faculty of  Math./Natural Sciences, and Interdisciplinary Centre for History and Philosophy of Science, \quad  scholz@math.uni-wuppertal.de} }
 \date{\small } 
 
\maketitle
\vspace{-4em}
\abstract{{\em This paper contains a case study of the work and self-definition of two important mathematicians  during the rise of modern mathematics: Felx Hausdorff (1868--1942)  and Hermann Weyl (1885--1955).  The two   had strongly diverging positions with regard to basic questions of mathematical methodology, which is  reflected in the style and content of their mathematical research. Herbert Mehrtens (1990) describes them as protagonists of what he sees as the two opposing camps of ``modernists''  (Hilbert, Hausdorff et al.) and ``countermodernists'' (Brouwer, Weyl et al.). There is no doubt that Hausdorff may  be described as a mathematical ``modernist''', while the qualification of Weyl as  ``countermodern'' is rather off the track, once his work is taken into account.}
}

\setcounter{tocdepth}{2} 
\tableofcontents         

\section{Introduction}
\label{sec: Intro}

In the history of recent mathematics there is a wide consensus that  mathematics  underwent a deep transformation in its epistemic structure and its social system in, roughly, the last third  of the 19th century and the first third of the following one. This  led to what is being considered as {\em  modern} mathematics in the sense of the 20th century. Jeremy Gray even called this  phase  a ``modernist transformation of mathematics''   \citep{Gray:Ghost}. His book  presents  a wide panorama of this period of change. The choice of the attribute ``modernist''  alludes to a wider cultural context of contemporary change in (visual, literary and sound) art and architecture. By good reasons Gray left it  open in what way, or even whether,  the   transformative tendencies in these different branches of culture can be comprehended as different expressions of a common historical phenomenon. This question is still wide open.\footnote{See the contributions by L. Corry and J. Ferreir\'os to this volume.}

By this choice of the word, and the question indicated by it, our author took up a suggestion of Herbert Mehrtens made in   {\em Moderne Sprache Mathematik} (sic!  without punctuation) \citep{{Mehrtens:Moderne}}. For me the title of this book  is difficult to translate, because it uses an ambiguity of the German language. It may be translated as ``Modern Language Mathematics'' or -- adding punctuation -- ``Modernity, Language, Mathematics''. Without the punctuation  the first alternative would be the correct translation, but one may also understand it in the second way.\footnote{The first reading  resonates with Mehrtens' way of presenting mathematics as a  language, organized in  two levels of {\em discourse}: the discourse {\em of} mathematics in the production, repectively documentation of knowledge, and the discourse {\em on} (about) mathematics, a meta-discourse which may include the foundational studies in the sense of Hilbert \citep[chap. 6]{Mehrtens:Moderne}. A non-anonymous referee (N. Schappacher) of the present paper opts for the second alternative.}
Mehrtens indicated the huge task of bringing together the historical understanding of the change in the  practice of mathematics as a social (and institutional) system (chap. 5)  and the knowledge style developed with it (chaps. 1--4). Apparently influenced by considerations from general history and history of art, he proposed to highlight  the radicality of the modern transformation of mathematics by establishing a  narrative of two opposing camps,\footnote{No social, political,  or cultural revolution ever happened without having to fight counter-revolutionary forces.}
 the protagonists of {\em modernity}  (``Moderne''), the  ``modernists'' driving  the modern transformation,  and  those opposing it, the ``countermodernists'', representing some (slightly mythical) entity  called {\em countermodernity} (``Gegenmoderne''). For both camps Mehrtens found two main, or at least typical, protagonists. David Hilbert (1862--1943) and Felix Hausdorff (1868--1942) (et al.)  for the modern camp versus Luitzen E.J. Brouwer (1881--1966) and Hermann Weyl (1885--1955) (et al.)  for the countermodern camp. Mehrtens made the  separation of the camps plausible by  arguing essentially on the  discourse level {\em about} mathematics, including the debate on foundational issues, with  only marginal references to mathematical knowledge (the mathematical discourse itself, to state it in his terminology). The strict separation of the camps  did not appear particularly  convincing to many readers and was {\em not} taken up by Jeremy Gray. But in a weakened form it seems to  remain a part of the  debate on modernism in mathematics.  The contributions of  L. Corry and J. Ferreir\'os to the present book  discuss this point from different perspectives; they also give good explanations of the terms  `` cultural modernism''  and ``modernism in mathematics'', which need not be doubled here.

 Hausdorff and Weyl, two protagonists of the opposing camps identified by Mehrtens, happen to  have been objects of my historical studies for some decades. The present book  is a good  occasion for laying down my view of  Mehrtens' presentation of these  mathematicians  as  representatives of his opposing camps.
 The following paper contains thus a simple case study  trying to  check  the adequacy and usefulness of the proposed categories for  our historical understanding of the 20th century.\footnote{For a discussion of Brouwer see the contribution of  Jos\'e Ferreir\'os  to this volume.} 
  Before drawing the conclusion in the last subsection of the paper, I try  to avoid, as far as possible, the qualifications  ``modernist'' or ``countermodern'', also of ``modernism'', at least as mathematics is concerned. On the other hand, I use the descriptive attribute ``modern'' for the deep transformation of mathematics between roughly 1860 and 1940 and the words ``rise of modernity'' for the social and cultural transformations in the late 19th century and the first half of the 20th. 
 
 The paper consists of three sections. It begins  with  informations on life and work of our protagonists in order to make the paper accessible for readers who are not so well acquainted with the details of the history of mathematics in the 20th century. Although the report on the mathematical work has to be extremely short at this occasion,  information on this aspect of our figures is of major importance for a judgement about the question how they stand in the process of the rise of  mathematical modernity. The section ends with an episode of parallel work of our authors on  Riemann surfaces in 1912. Both mathematicians had to deal with the problem  to make  Riemann surfaces more precise than was usual at  that time. The answers given to this challenge show  the different predilections  and different work styles of the two  in a nutshell (sec. 1.3). In the next section  the difference between our authors with regard to three questions which were important for both is discussed: the understanding of the mathematical continuum  (sec. 2.1), the relation between axiomatics, construction and the foundations of mathematics (sec. 2.2), and the question which role mathematics can, or ought to, play in the wider enterprise of understanding nature (sec. 2.3). The final section (sec. 3) discusses how our authors saw  themselves in the rise of modern society, before it comes back to the initial question of modernity and/or countermodernity.

\section{Hausdorff and Weyl, two representatives of 20th century  mathematics \label{section two mathematicians}}

\subsection{Two generations, two social backgrounds \label{subsection backgrounds}}
The main scientific work of our two protagonists took  place in the first half of the 20th century. Felix Hausdorff (1868--1942) was roughly  Hilbert's generation although seven years younger, Hermann Weyl (1885--1955) was a central figure of the next generation. 
Both were Germans, with Hausdorff coming  from a Jewish-German family.  During large parts of their life
both worked in German speaking countries.     Between 1913 and 1930 Weyl  lived in Switzerland,  
the liberal culture of which he learned to value against the torn and crisis stricken German social life during the inter-war years of the 20th century;  after 1933 he emigrated to the United States.   

Hausdorff obtained his doctorate (1891)   and habilitation  (1895)   at the University of Leipzig with   mathematical studies of the refraction and absorption of light as part of the research program of  the astronomer H. Bruns. He then turned towards  Cantor's theory of transfinite sets, the most abstract type of mathematics available at the time. Like many Jewish mathematicians of this time he  remained lecturer ({\em Privatdozent}) for a long time before he obtained his first associate  professorship ({\em Extraordinarius}) at Bonn University  in 1910. Three years later  he accepted a call to Greifswald  as full professor. In 1921 he returned to Bonn in the same position   and stayed there for the rest of his life.  Because of his Jewish  origin he  was in more immediate danger than Weyl after the rise of the Nazis to power, but  hesitated to emigrate in the early 1930s. When he finally tried so after 1939 he did not succeed. In 1942 when the anti-Jewish repression of the German Nazi regime reached its climax, he committed suicide with his wife and sister in law,  in order to elude deportation and death in the concentration  camps of the regime \citep{Purkert/Brieskorn:Hausdorff,Siegmund-Schultze:Hausdorffs-Tod}. 

Weyl had a different start into academic life. He obtained his  dissertation  (1908)  and habilitation (1910) in G\"ottingen   with research in real analysis (singular differential equations)  under the guidance of D. Hilbert and the intellectual influence of F. Klein. Swiftly accepted as a promising  young researcher he received a call as full  professor at the   Eidgen\"ossische Technische Hochschule Z\"urich already three years after his habilitation. In 1930 he hesitatingly accepted  a call to G\"ottingen as  successor of David Hilbert. In 1933, after the rise of the Nazi movement to power, he emigrated to the USA 
 following a call to the Institute of Advanced Studies where he was able to support other less privileged emigrants \citep{Siegmund-Schultze:Fleeing}. He stayed there until his retirement in 1951 and shuttled between Princeton and Z\"urich during the last years of his life. Only at rare occasion he visited post-war Germany. 

Both our protagonists came from well-to-do families and  grew up in German life and culture of the late 19th century and shared its humanistic higher school education.  Hausdorff's father was a successful textile merchant and owner of a small publishing house in Leipzig. As a traditional Jew he  participated in the community on a national level and was  active  against the  rising anti-semitism in late 19th century Germany. He stood in  opposition to the Jewish reformers  and contributed to the  scholarly orthodox  Talmud discussion  \citep[chap. 1]{Purkert/Brieskorn:Hausdorff}. Weyl's father was a director of a  local bank and city councillor of  Elmshorn, a medium sized town in Northern Germany. 
Already during his school time Weyl got deeply immersed in German philosophical  thinking by reading Kant's critique of pure reason in his parent home.\footnote{\citep[p. 632f.]{Weyl:EuB}}
 Although his first enthusiastic partisanship for a naive version of  Kantianism broke down  in the early years of his mathematical studies at G\"ottingen, in which he encountered Hilbert's axiomatic approach do geometry, he remained attracted by German idealistic philosophy in different molding, in particular Husserlian phenomenology and Fichtean constructive idealism \citep{Ryckman:Relativity,Sieroka:Umgebungen,Sieroka:neighbourhoods}. 

The young Hausdorff passed through a rather different intellectual trajectory. Standing in opposition to his father's orthodox Judaism he too was  attracted by Kant's critical philosophy, but he was also fascinated by Schopenhauer's pessimistic philosophy of life  and  of the young Nietzsche's radical cultural thoughts. In his later student's years he joined a circle of modernist intellectuals at Leipzig and participated and became active as a literary writer, essayist and free lance philosopher under the pseudonym Paul Mongr\'e. Different from Weyl, he considered  the liberation from any metaphysical bonds as a  desirable  goal of late 19th century thought \citep[chaps. 5, 6]{Purkert/Brieskorn:Hausdorff}. These differences in the general intellectual outlook between our protagonists would turn out to play a major  role for their predilections in mathematics and the way they reflected on their scientific work.

\subsection{Attempting the impossible: our authors' main contributions to mathematics
\label{subsection two contributors}}
Before we discuss the attitudes of our actors on methods, role and  goals of mathematical research  we have to recollect their main achievements in mathematics. We deal here with two ``giants'' of science and thus face an  essentially forlorn task, as  it would need book-length reports each to do justice to their work. Here we have to restrict to a selective, humble  survey on what may to be considered as the most important topics of their scientific achievements.  For  more extended report on F. Hausdorff see \citep{Purkert/Brieskorn:Hausdorff}, for Weyl \citep{Coleman/Korte:DMV,Chevalley/Weil:Weyl,Atiyah:Weyl}.

{\bf Felix Hausdorff} is well known for his axiomatization of the concept of topological space in his  opus magnum {\em Grundz\"uge der Mengenlehre} (Main Features of Set Theory) \citep{Hausdorff:Grundzuege}. But this book was much more.   Its first 6 chapters contained the leading introduction to Cantorian set theory in the first decades of the 20th century and included a detailed study of transfinite order structures to which Hausdorff had contributed himself in the preceding years  (including the study of $\eta_{\alpha}$ sets which later became important in  foundational studies of set theory). The second half of the book   
established a program for founding  basic fields of mathematics on an axiomatization in the framework of set theory. As we now understand,  this was a main trend of the modernization of mathematics in the first half of the 20th century culminating in the work of ``modern algebra''  and  Bourbaki's vision of mathematics around the middle of the century. Hausdorff himself exemplified the method for topological spaces (chap.7), metrical spaces (chap. 8), functions (chap. 9),   measure theory and integration (chap. 10). In these chapters he could in particular draw the consequences of developments of the first abstract (topological) space concepts  about the turn of the centuries, due to M. Fr\'echet, F. Riesz, E.H. Moore and others.\footnote{See  \citep[p. 353ff.]{Purkert/Brieskorn:Hausdorff};  a more  extended discussion of the early development of topological space concepts is given in the commentary (in German) on the historical background of Hausdorff's axioms  in \cite[vol. 2, pp. 675--708]{Hausdorff:Werke}.}

At the end of the book Hausdorff published a  paradoxical disjoint decomposition of the 2-sphere (using the axiom of choice), republished separately in  \citep{Hausdorff:Paradox}. 
According to his own description it showed that in Euclidean space ``one third of the sphere''  may be congruent to   ``one half'' (ibid, p. 430).
Hausdorff was a master of logical precise argumentation without using a formal system for logic itself  and was fond of counter-intuitive effects in the world of transfinite set theory.    Some years later he generalized a measure theoretic approach initiated by C. Caratheodory and introduced a class of measures on subsets of metric spaces ({\em Hausdorff measures}) which allow to characterize the dimension $p$  of point sets, where $p$ may assume {\em fractional} values \citep{Hausdorff:Masse}.  
 The concept developed in this small paper on measure and dimension has been of enormous influence in mathematics (non-linear partial differential equations, dynamical systems, ergodic theory) and physics (potential theory, turbulent flows, meteorology) and has become widely known through with the rise of fractals 
and computer graphics  at the end of the 20th century 
 
In his investigations of set theory he introduced important fundamental concepts like co-finality and  co-initiality of ordered sets, or (Hausdorff) gaps in dense ordered structures. After  establishing the  distinction of regular and singular initial (cardinal) numbers  he
observed  that regular cardinal numbers with limit index, should they  exist at all,   would be of ``exorbitant'' size  \citep[233]{Plotkin:Hausdorff}. Later this became the starting point for the study of so-called ``large cardinal numbers''.  Among his diverse  contributions to set theory we also find the Hausdorff maximality principle of partially ordered sets, the  general recursion formula of aleph-exponentiation, and the concept of so-called $\eta_{\alpha}$ sets which later became important in model theory.\footnote{This list  is a  
 a selection of the survey of Hausdorff's main contributions to set theory in  \citep[p. x]{Purkert/Brieskorn:Hausdorff}. For more details see there.}

Other contributions of his relate to different fields of mathematics of the 20th century, e.g., the Baker-Campbell-Hausdorff formula in the theory of Lie groups, or the Hausdorff-distance of compact subsets of a  metrical space, which was later used by M. Gromov to measure the ``distance'' of metrical spaces from being isometric. The resulting Gromov-Hausdorff distance of metrical spaces became an important tool for differential topology  \citep[vol. 1B, p. 779]{Hausdorff:Werke}.   In Hausdorff's lecture manuscripts we find many interesting insights, e.g. with regard to an axiomatic foundation  of probability theory \citep[vol. {\em 5}, 595--723]{Hausdorff:Werke} similar to the one in Kolmogorov's famous book of 1933, \nocite{Kolmogoroff:1933} but   ten years earlier. The mentioned topics show already the  profile of Hausdorff's contribution to 20th century mathematics: set theory as the basis for work and as a framework of modern mathematics, order structures, point set topology, metric spaces, measure theory with particular attention to paradoxical or seemingly paradoxical (fractional dimension) results, and functional analysis.   After his turn towards pure mathematics the  contributions to applied mathematics of his early phase were no longer of interest to him. On the other hand,  his interests in pure mathematics show  traces of his epistemological  reflections  around 1900  in which he assigned mathematical arguments an important role for  the critique and decomposition of classical metaphysics \citep[sec. 5f.]{Epple:Hausdorff2021}.  Hausdorff's interest and active  participation in  philosophical reflection of  mathematics faded away after his turn towards pure mathematics research. As  Purkert/Brieskorn write, his alter ego Paul Mongr\'e ``bid farewell to the public'' about 1910 \citep[p. 318]{Purkert/Brieskorn:Hausdorff}.

{\bf Hermann Weyl}, on the other side, was acknowledged as a leading figure of the post-Hilbert generation of mathematicians already during his life time. His research was as broad as the one  of his academic teacher Hilbert; it comprised many fields inside mathematics and its foundations as well as long lasting contributions to mathematical physics. He  was widely read in philosophy and did not hesitate to  share his philosophical reflections on mathematics and science with the interested public. His most influential work in mathematics proper results from  his studies in Lie theory  starting in the mid-1920s \citep{Weyl:Darstellung}. He  combined E. Cartan's characterization of infinitesimal groups (Lie algebras) and their representations with an integral approach  used by  I. Schur to the characters of certain  groups (the special orthogonal ones). Weyl was able to generalize Schur's method to  all the classical groups and  to study their representations \citep{Hawkins:LieGroups,Eckes:Diss}. In his Princeton years he extended this approach, in cooperation with R. Brauer, to give a modern access to the invariants of the classical groups \citep{Weyl:ClassGroups}.

Although lying deep inside mathematics proper (i.e., ``pure''  mathematics),  for Weyl this research topic was multiply  intertwined with questions coming from theoretical physics and leading back to the latter. He had started to develop interest in infinitesimal symmetries in his thoughts about general relativity and  generalized Riemannian geometry by introducing what he called a ``length'' gauge (today scale gauge). This led him to propose a  geometrically unified field theory of gravity and electromagnetism in the framework of the first gauge theory of electrodynamics  with local symmetries of geometric scale  as the {\em gauge group}  \citep{Weyl:GuE,Weyl:InfGeo}.\footnote{In Jed Buchwald's contribution to this volume   ``gauging'' is discussed in the  pre-Weylian perspective of under-determination  of the  electromagnetic potential (decomposed in its scalar and its vector part) up to exact differentials as  a history of ``gauge'' {\em ante letteram}. We learn from it that important physical questions of this under-determination have been posed and answered long before the explicit concept of {\em gauge} was introduced; for the later development see, among others,  \citep{ORaifeartaigh:Dawning,ORaif/Straumann:2000}. } 
In this context he  contributed importantly  to  clarifying  conceptual and mathematical questions in general relativity \citep{Weyl:RZM} and  differential geometry \citep{Scholz:Connections,Scholz:ICM}.

His proposal of a scale gauge theory of electromagnetism did not work out directly as a physical theory  but  could be ``recycled'' after the advent of the new quantum mechanics in  form of a gauge theory for the phase of wave functions of charged particles \citep{Vizgin:UFT,Scholz:Weyl_PoS}. Through the intermediation of W. \cite{Pauli:1933}, \citep{Pauli:1941}, Weyl's idea of a  gauge field approach to electromagnetism was generalized by C.N. Yang and L. Mills in 1954 to a more general gauge group of isotopic spin $SU(2)$ \citep{ORaif/Straumann:2000}. After a long interlude of laborious research in high energy  physics it acquired a  central role in the standard model of elementary particle physics in the 1970 \citep{ORaifeartaigh:Dawning}. About the same time it  entered also the research of differential topology and was used for defining new topological invariants \citep{Kreck:1986}. 
On a different, although connected route Weyl started to use group representation theory in the new quantum mechanics after 1926 \citep{Weyl:GQM}. Together with E. Wigner he may be considered as a main actor for propagating  symmetry considerations in the study of quantum systems which again became a main tool for particle physics in the second half of the 20th century \citep{Borrelli:isospin,Borrelli:2015}. 

A third field in which Weyl intervened with long lasting consequences was  the debate on the foundations of mathematics  in the first third of the 20th century. In spite of his high regard for Hilbert as a mathematician he was not at all convinced philosophically by his teachers proposal for a formalistic solution of the problems arising in transfinite set theory  around 1905 and the consequences for analysis, arithmetic and mathematics in general. Weyl started to develop a constructive alternative for the foundation of analysis \citep{Weyl:Kontinuum}. A little later he   even fought for some years at the side of Brouwer for an  intuitionistic program  in the foundations of mathematics attacking Hilbert harshly \citep{Weyl:Krise}. In later writings he came to a more balanced view of Hilbert's foundational program (see below).

Like in the presentation of Hausdorff's work this survey is necessarily extremely selective: other fields of Weyl's work, e.g.,  in convex geometry, real and complex  analysis  have fallen completely  through the cracks. The next subsections gives the chance for partially correcting  this at least with regard to complex analysis.

\subsection{Contrasting trajectories: Riemann surfaces as an example \label{subsection Riemann surfaces}}
By a funny historical coincidence  our two protagonists lectured on complex function theory at the same time without knowing of each other. In the  winter semester 1911/12, Hausdorff gave an introduction to function theory  at Bonn university, Weyl  lectured on Riemann's theory of Abelian integrals in  G\"ottingen. Both had to struggle with the  concept of a Riemann surface, which at that time was still only vaguely defined,  and both made proposals how to attack this question, with  long ranging repercussions.  Weyl's notes became a draft for his book  on the idea of  Riemann surface  (``Die Idee der Riemannschen Fl\"ache'') published in the following year \citep{Weyl:IdeeRF}. This book is widely known for presenting the first definition of a  manifold at least for the 2-dimensional case. For  Hausdorff the lecture gave him reason to think about neighbourhood systems which turned into his axiomatics of topological spaces two years later.

 Weyl drew upon Hilbert's sketch of an axiomatic characterization of the (real) plane, based on  topological concepts \citep{Hilbert:1903Grundlagen}.\footnote{Also in \citep[appendix IV]{Hilbert:GG2}.}
Hilbert defined a {\em plane} as a ``system of things''  (set), with elements (``things'') called {\em points}, which is bijectively mapped as a whole  on the ``number plane''. He then used Jordan domains of the latter for characterizing {\em neighbourhoods} (``Umgebungen'') of points in the plane. Weyl could link up to this idea but had to modify it.  For Riemann surfaces, thought to arise semi-constructively from analytic function elements (``analytische Gebilde'') in the sense of Weierstrass, he  had to localize Hilbert's idea and could no longer presuppose a global bijection with the number plane. This led  to the first definition of a manifold $\mathfrak{F}$ in Riemann's sense, although restricted to the 2-dimensional case, by establishing an axiomatics of neighbourhood systems  in  $\mathfrak{F}$, with bijective maps to  open disks in  the Euclidean plane \citep[p. 17f.]{Weyl:IdeeRF}. This sufficed for  defining  continuity,  differentiability and even analyticity of maps between such manifolds and of functions and to build Riemann's theory of Abelian differentials on such a fundament. 
In particular  the topological  notions of triangulation, simple connectedness,  covering surfaces, group of covering transformations and the topological genus of the surfaces etc. were thus put on an essentially clarified mathematical basis, if one kept the constructive context (analytic function elements and disks as coordinate images) in mind. A later analysis from the more refined point of view of Hausdorff's topology would show that Weyl's axioms were not strong enough as a self-sufficient formal axiomatization: the later Hausdorff separation property was not secured by his axioms although implicitly presupposed in the derivations. But this was not Weyl's concern  during the next few decades. Only during the preparation of the English translation by lectures given in 1954 at Harvard and Princeton, and in  the third German edition Weyl finally added  Hausdorff separation as an further axiom \citep[p. xii]{Remmert:WeylRF}. 

 The idea to talk about {\em neigbourhoods} of points  not only in geometry proper but also in Weierstrassian function theory and even  for characterizing more general spaces on the background of set theory was not an exclusive privilege of the G\"ottingen mathematicians. Weierstrass had used the terminology already, and also F. Riesz used it in thoughts about generalized spaces \citep{Riesz:1908,Rodriguez:Diss}. Hausdorff, who  had started to lecture on Cantorian set theory in Leipzig in summer semester 1901 and again in Bonn in 1910, had not yet taken up this idea  in his discussions of topological aspects in general sets. In his lectures on function theory he was confronted with neighbourhoods of function elements in a natural way. In his lecture notes of winter semester 1911/12 we find clear evidence that  he realized at this point that the study of neighbourhood systems was the clue for ``ordering the system of points''  which arose  from the study of equivalence classes of analytic function elements, and also more generally. Moreover he became aware that the structural properties of such neighbourhood systems had to be analysed.  This he did more extensively  in summer semester 1912, in which he gave his next course on set theory at Bonn. Here we find four structural properties of systems of neighbourhoods in metrical spaces which were essentially the axioms of topological space, published in his book \citep{Hausdorff:Grundzuege}. Moreover he already announced that these structural properties could be used as axioms for general spaces  \citep[p. 714ff.]{Epple_ea:TopRaum}.  
 
 This small episode seems characteristic for the different thought styles of our protagonists:  Hausdorff used the analysis of  conceptual features of Weierstrassian function elements underlying the concept of  Riemann surface  as a stepping stone for a more fundamental search of a  general  characterization of topological spaces in the framework of  sets.  Weyl, on the other hand,  took the same incentive as a starting point for establishing  an axiomatic clarification of the  intuitive concept of Riemann surface which had been  in use already for several decades. He kept it   closely linked to the construction of global objects from Weierstrassian function elements, aiming at concrete mathematical objects with multiple structures. This difference may seem a nuance of research orientation only; but we will see that it is characteristic for their contrasting positions with regard to the aim and character of mathematics.  During the following years it  would  develop  into an open opposition.

\section{Mathematics in the tension between formal thought and   insight  \label{section formalism etc}}

\subsection{Two opposing views of the  {\em continuum}:   a modified classical concept versus a set theoretic perspective \label{subsection set theory}}
Already in the period 1910 to 1914 Hausdorff and Weyl had developed quite different ideas of how to deal with the mathematical concept of {\em continuum}. As we know already Hausdorff was attracted by the epistemic perspective opened up by Cantor's transfinite set theory once he got to know of it; in contrast to this  Weyl became increasingly sceptical with regard to any truth claim for the latter  after discussing  foundational problems in his Habilitation lecture (1910) (see below).  W. Purkert was able to reconstruct from indirect evidence (remarks on the infinite in philosophical essays written under the pseudonym Mongr\'e) that Hausdorff got to know Cantor's theory during the year 1897, the year of the First  International Congress of Mathematicians at Z\"urich \citep[p. 262ff.]{Purkert/Brieskorn:Hausdorff}. Hausdorff/Mongr\'e was fascinated by the intellectual perspective of Cantor's treatment of the transfinite cardinal and ordinal numbers (although not yet clarified in sufficient detail, not to speak of its axiomatization). At this time he pursued a  philosophical program in the footsteps of Kant,  radicalized by Nietzsche, for ``proving'', more precisely by arguing with the use of mathematical metaphors,  that no knowledge of the  ``thing in itself'' is possible, and in particular no insight into the structure of  ``absolute time'' or  ``absolute space'' or even ``cause''  is possible \citep{Stegmaier:Hausdorff,Epple:Hausdorff2021,Hausdorff:Chaos,Hausdorff:ChaosSelbstanzeige}.  Wearing the hat of Mongr\'e,  our author tried to convince his readers by an ``apagogic proof'' (a proof by contradiction)  that  absolute time or space, if assumed,  cannot have any type of structure. For this goal Cantorian set theory seemed to him an ideal tool. Relative structures, i.e., not completely absolute ones,  were of course possible also for him, i.e., order structures in the case of time and geometrical or even topological structures (before the advent of the word) in the case of space. According to Mongr\'e/Hausdorff such non-absolute structures  were ``selected'' by the mind, to make human action and survival possible; they could then form a  rather individualistic ``cosmos''. The individualistic  exaggerations to be found at many places of the early  Mongr\'e's literary and some of the philosophical utterances were, however, step by step  moderated and substituted  by what Hausdorff  a bit later called a {\em considerate empiricism}, which respected  empirically founded scientific knowledge, including theoretical refinement and critique  \citep{Epple:Hausdorff2006}. 

As M. Epple and other authors have argued, Hausdorff's mathematical research topics  in the time between 1900 and 1914 was still embossed by his interest in logically consistent, although intuitively surprising, perhaps even paradoxical insights in order structures ($\sim$ time) and/or topological,  metrical and measure structures ($\sim$ space) \citep[sec. 5.6]{Epple:Hausdorff2006}. Hausdorff's first important works in transfinite set theory consisted in profound and technically demanding contributions to order structures \citep{Hausdorff:Ordnungstypen}. W. Purkert observed that an additional motivation for this work seems to have been to come closer to a proof of Cantor's famous {\em continuum hypothesis}; i.e. the assumption (at first a claim of Cantor) that the cardinal number of the subsets of the natural numbers $2^{\aleph_0} = \mathfrak{c}$  (which encodes the cardinal number of the ``continuum'', i.e. the real numbers)  is the first non-denumerable cardinal number 
\[ 2^{\aleph_0} = \aleph_1
\]
or, even more generally,  $2^{\aleph_{\nu}} = \aleph_{\nu+1}$.

In Hausdorff's view, the ``continuum'' itself would have to be understood by using  all kinds of different types of topological and/or measure structures. Intuitive insight into the nature of the continuum seemed him of ephemeral value only, perhaps important for the imagination of the individual mathematician, but without any epistemic value with regard to truth claims. His great book {\em Grundz\"uge der Mengenlehre} was a splendid exemplification of this general view.

Weyl had a completely different view of the continuum,  which was deeply influenced by  the long tradition in mathematical and philosophical thought upon this subject. Riemann's concept of {\em manifold} appeared to him as the most promising modern clue to the topic. Its logical and formal foundations remained an open question for him until the end of his life, although he himself  made at least three attempts to come to grips with it  \citep{Scholz:2000Continuum}: a constructive approach  in \citep{Weyl:Kontinuum}, influenced by E. Borel and H: Poincar\'e, 
which was designed to avoid the pitfall of impredicative definitions,\footnote{For an appraisal of its mathematical long range import see \citep{Feferman:Weyl_Kontinuum}; a critical historical view is given in \citep{Schappacher:points}.}  an intuitionistic one in \citep{Weyl:Krise}, and a combinatorial topological one at different occasions \citep{Weyl:Analysis_situs} or in his lecture course on {\em Axiomatics} in G\"ottingen 1930/31 \citep[\S 37]{Weyl:Axiomatik1930}. 

In a paper written for the Lobachevsky anniversary in 1925, though published only posthumously, we find a most explicit remark why Weyl would not agree with Cantor's or Hausdorff's approach to the continuum, at least  understood in the sense of a manifold describing physical space. In such a manifold the local descriptions by coordinates in a ``number space'' are ``arbitrarily projected into the world'' and everything else, in particular the metric structure is turned into a field on the space. This could  well be reflected in a Kantian type of approach which shaped  his understanding of the role of the spacetime concept in general relativity.  A few years earlier he had contributed to a deeper understanding of the underlying concepts by posing  
the problem of space anew, facing the changed situation after the rise of special and general relativity \citep{Scholz:2016Weyl/Cartan,Bernard/Lobo:PoS}. In the mid-1920s Weyl  resumed a (relativized) Kantian perspective and sharpened his criticism of a set theoretic substitute for it in the following way:
\begin{quote}
Space thus  emerges  [by separating the topological manifold from the metrical and other fields on it, ES] even more clearly as the form of appearances in contrast to its real content: the content is measured once the form has been referred to coordinates. [Set theory, one may say, goes even further; it reduces the mf [manifold] to a set as such and considers already the continuous connection  as a field on the latter. It should, however, be clear that in doing so  it violates against the essence of the continuum which by its nature cannot be smashed into  a set of isolated elements. The analysis of the continuum should not be founded on   the relation between element to the set, but on the one between part and the whole.  \ldots] \citep[p. 4f., second square brackets in  orig., translation ES]{Weyl:RiemannsIdeen}\footnote{``Deutlicher tritt dadurch der Raum als Form der Erscheinungen seinem realen Inhalt gegen\"uber:
der Inhalt wird gemessen, nachdem die Form willk\"urlich auf Koordinaten bezogen ist. [Die Mengenlehre, kann man sagen, geht darin
noch weiter; sie reduziert die Mf auf eine Menge schlechthin und
betrachtet auch den stetigen Zusammenhang schon als ein in ihr
bestehendes Feld. Es ist aber wohl sicher, da\ss{}  sie dadurch gegen das
Wesen des Kontinuums verst\"o\ss{}t, als welches seiner Natur nach gar
nicht in eine Menge einzelner Elemente zerschlagen werden kann.
Nicht das Verh\"altnis von Element zur Menge, sondern dasjenige des
Teiles zum Ganzen sollte der Analyse des Kontinuums zugrunde
gelegt werden. Wir kommen darauf sogleich zur\"uck.]'' \citep[\S 37]{Weyl:RiemannsIdeen}}
\end{quote}  
With other words, Weyl considered transfinite sets  as an overstretched formal concept without substantial content, at least as far as physical spacetime is concerned.\footnote{Ferreir\'os (2016)\nocite{Ferreiros:2016practices} calls this a ``pointillist'' view of the continuum, see in particular the discussion in chap. 10.4.} Below we see that his scepticism did not only relate  to the continuum as a concept of mathematical physics but also  to its role  in   mathematical analysis and  in the foundations of mathematics.

\subsection{Axiomatics,  construction and the open problem of the  foundations of math \label{subsection axiomatics}}
 Weyl understood axiomatics as  the defining basis of a conceptual framework on which a mathematical theory could be built. He saw no opposition between axiomatics and the construction of mathematical object fields. The task of an axiomatics formulation was  to clarify the structure of some field  of mathematical thought; its objects were to be constructed and  dealt with  symbolically. This should happen  without too strong hypotheses about the infinite, in particular without making use of transfinite set theory and only if unavoidable with applying the principle of the excluded third without a constructive underpinning. Weyl's axiomatization of the 2-dimensional manifold and of Riemann surfaces was an early example. And he stuck to this conception essentially for his whole life (i.e., with gradual modifications only). In the  late 1930s  he came into contact with members of the early Bourbaki group, in particular C. Chevalley, and started to develop more respect for the algebraists usage of axiomatics  as a research tool in its own right, although still not  with respect to foundational issues (which were not in the focus of Bourbaki anyhow).
 
  Weyl was extremely sceptical with regard to Hilbert's program of founding mathematics by axiomatization and a  formal analysis of the proof structure with the aim of showing its internal consistency. He considered such a justification of mathematical theories or even of mathematics as a whole   as nothing but a  formalistic showpiece   which might be impressive because of its acumen,  but  would fail completely the goal of justifying the substance of mathematics. 
In his view a meaningful justification would presuppose a clarification of the basic conceptual ingredients of a mathematical theory by  {\em symbolic construction}, as he called it  \citep{Weyl:PMN,Weyl:PMNEnglish}. 
A preliminary version of how a constructive  approach to analysis might work was given in   his famous book {\em Das Kontinuum} \citep{Weyl:Kontinuum}; for its long ranging impact see \citep{Feferman:Weyl_Kontinuum}. But Weyl was discontent with his own achievements, not justbecause it justified only a restricted variant of analysis (without the general existence of a supremum of a bounded set of the reals). After he had constructed his reduced (denumerable) range of real numbers, he  opened the discussion of the relation to geometry with a self-critical remark. 
He deplored that the intuition of connectivity inherent in the  geometrical concept  {\em continuum} was not depicted in his constructive number continuum: 
\begin{quote}
Once we have torn the continuum apart into isolated points, it
is difficult to  reconstruct  ex post 
 the connectivity between the single points, which is based on their non-independence,    by  some conceptual equivalent. \citep[79, translation ES]{Weyl:Kontinuum}\footnote{``Nachdem wir das Kontinuum in isolierte Punkte zerrissen haben, f\"allt es jetzt schwer, den auf der Unselbst\"andigkeit der einzelnen Punkte beruhenden Zusammenhang nachtr\"aglich durch ein begriffliches \"Aquivalent wieder herzustellen'' \citep[79]{Weyl:Kontinuum}. The translation in \citep[103f.]{Weyl:KontinuumEnglish} suppresses the details ``ex post'' and the proxy character of the ``conceptual equivalent''.}
\end{quote}
So his constructive (denumerable) continuum of 1918 offended against the ``essence'' of the continuum at least as much as a Hausdorffian set theoretic approach (criticized in the quotation at the end the last subsection). Its only advantage was  its (semi-finitist) constructive methodology rather than the  one in which transfinite sets were postulated axiomatically. 
Irrevocable connectivity between points by  their inseparable infinitesimal neighbourhoods was what Weyl looked for. For a while he believed to find it in the intuitionist approach proclaimed by L.E.J. Brouwer more or less at the same time \citep{Brouwer:1919}. So Weyl's attempts at  laying the cornerstones of  a constructive clarification for analysis shifted for some years (between 1919 and 1923) towards a strong  support for  Brouwer's more radical intuitionistic program, most decidedly expressed in his open polemics of \citep{Weyl:Krise}. This most radical phase of his contributions to the foundations of mathematics has attracted much attention in the history and philosophy of mathematics  \citep{Rowe:Weyl,Rowe:BrouwerHausdorff,Hesseling:Gnomes,%
Scholz:2000Continuum},  \citep[sec. 4.1]{Mehrtens:Moderne}, and with more technical details \citep[sec. 6]{Coleman/Korte:DMV}. 
 
About the mid-1920s  he started to  accept that Hilbert's formalist program  was, after all,  a defensible position. He remained sceptical, however, with regard to the epistemic value of such a  formalist axiomatic approach to the foundations of mathematics and to the concept of continuum, because it did not live up to his (undefined and probably  undefinable) criteria  of ``insight'' and ``meaning''. 
  At different occasions he sketched how he would imagine a  constructive symbolic  approach to the continuum,  based on  methods taken from combinatorial topology. He explored here how far a symbolic representation of cell complexes with (denumerably) infinite sequences of barycentric subdivisions would carry  \citep{Weyl:Analysis_situs,Weyl:Axiomatik1930,Weyl:math-thinking,Weyl:1985}; but a purely combinatorial constructive characterizations of topological manifolds, which he considered as the best mathematical approach to the ``continuum'',  remained an unsolved problem. 
  
   At the turn towards the middle of the century he accepted and appreciated the meanwhile widely spread axiomatic approach in mathematics: 
  \begin{quote}
  \ldots the axiomatic attitude has ceased to  be the pet subject of the methodologists [researchers in formal logic and foundations of mathematics, ES] its influence has spread from the roots to all branches of the mathematical tree \citep{Weyl:math-thinking}.
    \end{quote}
But it remained important for him that axiomatic postulates were not dissolved from ``symbolic construction''. He was not satisfied with taking   finitist methodology  serious {\em only} at  the level of metatheorical investigation  (like Hilbert had proposed in his  proof theoretic program for showing the consistency of axiomatic theories). He demanded that on  {\em  all  levels}  of its knowledge production and  reflection   mathematics ought to be a 
``\ldots dexterous blending of constructive and axiomatic procedures''  \citep[38]{Weyl:1985}.  The foundations of mathematics, on the other hand,  remained an open problem for him until the very end of his life. \\[0.2em]

In one respect Hausdorff's view of axiomatics was  not too different from Weyl's (and any other 20th century mathematician): in the sense of giving the format for defining the basic concepts of a mathematical theory. But in others it differed drastically.  Hausdorff was an excellent logically sharp thinker who did not see a  need for  formalizing logic, as we noticed already. He would never give up a  principle like the one of the excluded middle. This would unnecessarily  reduce the range of mathematics  and was completely unacceptable to him.  Symbolically supported  {\em creation} combined with logical precision took the place  occupied by symbolic {\em construction} for Weyl in the generation of mathematical knowledge. Transfinite set theory as outlined by Dedekind and Cantor, continued by himself and adopted by Hilbert and his school 
 established a language and thought milieu for symbolic creativity par excellence. He knew, of course, about the open questions in the foundations of set theory, in particular  that the comprehension of infinite totalities had to be handled with care, but he saw no hindrance to build mathematical theories along these lines. Still in the late 1920, at the occasion of the so-called second edition of his book on set theory (in fact, a new book) he emphasized the aspect of creativity in  his lucid rhetoric:
 \begin{quote}
 It is the  eternal achievement of Georg Cantor to have dared this step into infinity, under interior and exterior struggles against seeming paradoxies, popular prejudices, philosophical statements of power (infinitum actu non datur), but also against reservations pronounced by the greatest mathematicians. By this he has become the creator of a new science, set theory, which today forms the grounding of the whole of mathematics. In our opinion, this triumph of Cantorian ideas is not belittled by the fact that a certain antinomy arising from an excessively boundless freedom of forming sets still needs a complete elucidation and elimination. \citep[11, {Werke {\em 3}, 55}]{Hausdorff:Mengenlehre2}\footnote{``Es ist das unsterbliche Verdienst Georg Cantors, diesen Schritt in die Unendlichkeit gewagt zu haben, unter inneren wie \"au\ss{}eren K\"ampfen gegen scheinbare Paradoxien, popul\"are Vorurteile,
philosophische Machtspr\"uche (infinitum actu non datur), aber auch gegen Bedenken, die selbst von den gr\"o\ss{}ten Mathematikern ausgesprochen waren. Er ist dadurch der Sch\"opfer einer neuen Wissenschaft, der Mengenlehre geworden, die heute das Fundament der gesamten Mathematik bildet. An diesem
Triumph der Cantorschen Ideen \"andert es nach unserer Ansicht nichts, da\ss{} noch eine bei allzu uferloser Freiheit der Mengenbildung auftretende Antinomie der vollst\"andigen Aufkl\"arung und Beseitigung bedarf.'' \citep[11, {Werke {\em 3}, 55}]{Hausdorff:Mengenlehre2}
}
 \end{quote}
 
He knew that in the environment of Hilbert (Zermelo, Fraenkel, Bernays) the axiomatization of set theory was under way, but saw no pressure to              proceed along these lines, and was far from feeling any ``anxiety'' that something would go wrong with the foundations set theory and mathematics \citep[sec. 7.3]{Purkert/Brieskorn:Hausdorff}.

From such a viewpoint a methodology which demanded a reduction of symbolic creation to procedures that would deliver only denumerable ranges of objects (Weyl's constructivism or Brouwer's intuitionism) appeared to him ridiculous, or even worse. He did not state such an opinion publicly, but was clear up to extremity in a letter to Abraham Fraenkel written  June 9, 1924, in  response to Fraenkel's step forward with regard to the foundations of set theory   \citep{Fraenkel:1923}. He thanked for the progress his correspondent had achieved for an axiomatic framing of set theory and the discussion of the set theoretic antinomies, because this spared him the work to deal with  questions for which he has no knack for ("Dinge, die mir nicht liegen''). From now on he would be able just to refer to Fraenkel's book. He  continued:
\begin{quote}
You have even succeeded in making the oracle pronouncements of
Brouwer and Weyl understandable -- without making them appear to me
any less nonsensical! You and Hilbert both treat intuitionism with too
much respect; one must for once bring out heavier weapons against the
senseless destructive anger of these mathematical Bolsheviks! \citep[vol. {\em 9}, 293, translation \citep{Rowe:BrouwerHausdorff}]{Hausdorff:Werke}\footnote{``Es ist Ihnen sogar gegl\"uckt, die Orakelspr\"uche der Herren Brouwer und Weyl verst\"andlich zu machen -- ohne dass sie mir nun weniger unsinnig ercheinen! Sowohl Sie als  auch Hilbert behandeln den Intuitionismus zu achtungsvoll; man m\"usste gegen die sinnlose Zerst\"orungwuth dieser mathematischen Bolschewisten einmal gr\"oberes Gesch\"utz auffahren! \ldots'' \citep[vol. {\em 9}, 293]{Hausdorff:Werke}}
\end{quote}
Hausdorff's surprisingly militant language has to be understood on the background of Weyl's polemical language  in  his paper propagating an intuitionistic ``revolution'' -- and the excited time conditions in post-war Germany of the early 1920s. It indicates a deep dividing line among early 20th century mathematicians (in Germany) with regard to basic methodological convictions and the value of certain research programs. But can this dividing line be better understood by declaring our two protagonists as belonging to two separate camps of {\em modernists} (Hausdorff) and {\em counter-} or even {\em antimodernists} (Weyl) as proposed by \citep{Mehrtens:Moderne}? -- David Rowe calls Weyl, just to the contrary,  a ``reluctant revolutionary'' \citep{Rowe:Weyl}. This seems to me much more to the point; we will come back to this question in the final discussion. 

\subsection{Mathematics and the  material world  \label{subsection applied math`}}
  Although Hausdorff did no longer contribute actively to natural sciences after his disappointing experiences with his early works in astronomical optics, he  held a pronounced opinion with regard to the question which role mathematics may play for understanding the outer world via its use in the natural sciences  \citep[chap. 3]{Purkert/Brieskorn:Hausdorff}. His  contributions to probability theory remained relatively unnoticed  \citep[sec. 4.2, 10.3 ]{Purkert/Brieskorn:Hausdorff}; the lecture containing a set theoretic axiomatization of probability  remained unpublished \citep{Hausdorff:VorlesungWT}.

 In the 1890s and early 1900s he was highly interested in the question of non-Euclidean geometry and in philosophico-mathematical question of space and time concepts \citep{Epple:Hausdorff2021}, \citep[sec. 5.6]{Purkert/Brieskorn:Hausdorff}. In his radical thoughts on philosophical (epistemological and ontological) questions, published under the name Paul Mongr\'e, mathematics played an important role for undermining the belief in fixed, perhaps even a priori, forms of knowledge of the external (material) world. The great variety of geometrical or,  in nuce, even topological structures for spacelike thinking, and of order structures for timelike thinking became an important tool for him in putting established notions of mathematical physics, astronomy and cosmology in question. On the other hand, he made sure that the ordering of sense perceptions and scientific empirical knowledge needed mathematics for acquiring a well defined and intelligible form. He called such a methodology  {\em considered empiricism} (``besonnener Empirismus''), in contrast to empiricism sans phrase and positivism on the one hand and neo-Kantianism, or any other rejuvenated version of German idealism, on the other \citep{Epple:Hausdorff2006}. 

In later years Hausdorff did not completely  lose  interest in mathematical physics but it clearly moved to the background of his attention. We know that he prepared talks, perhaps even a introductory publication for a wider public, on (special) relativity \citep[sec. 5.1]{Epple:Hausdorff2021}, but he never took up questions from mathematical physics for his own research. From the beginning of the 20th century onward his research profile became the one of a ``pure'' mathematician who appreciated the role of mathematics for an  open minded and critical understanding of the material world. In his early years he had formulated a basic attitude underlying such a role of his science in an aphorism:
\begin{quote}
What we are  missing is a self-critique of science; the verdicts of  science given by   art,  religion and  sentiment  are just as numerous as useless. Perhaps this is the ultimate destination of {\em mathematics}. \citep[aphorism 401, transl. ES]{Mongre/Hausdorff:St-Ilario}  \footnote{``Uns fehlt eine Selbstkritik der Wissenschaft; Urtheile der Kunst, der Religion, des Gef\"uhls \"uber die Wissenschaft sind so zahlreich wie unn\"utz. Vielleicht ist dies die letzte Bestimmung der {\em Mathematik}.''   \citep[Aph.  401]{Mongre/Hausdorff:St-Ilario} }
\end{quote}

Weyl, on the other hand, was a highly creative contributor to mathematical and theoretical physics, besides his great achievements in pure mathematics, by far too huge to be resumed here. As is well known, he  made outstanding contributions to Einstein's theory of gravity and early cosmology  \citep{Giulini:RZM2023,Lehmkuhl:Einstein-Weyl,Goenner:DMV,Rowe:Cosmology-II}, the generalization of Riemannian geometry as a scale covariant (conformal) framework for relativistic field theory \citep{Vizgin:UFT,Ryckman:Relativity,MacCoy:2021}, to the introduction of the gauge principle into the rising quantum mechanics \citep{Straumann:Einstein_Weyl,ORaifeartaigh:Dawning},
 and  finally he  displayed, conjointly with B.L van der Waerden and E. Wigner,  the usability of group representations as a basic frame for studying symmetries in quantum physics \citep{Eckes:Weyl_Wigner,Schneider:Diss,Scholz:WeylGQM}. 
All of this  turned out to be of long ranging influence on the course of physics during the 20th century, and probably also beyond  \citep{Yang:Weyl,Mackey:WeylsProgram,Borrelli:2015,Borrelli:selectionrules,%
Scholz:2018Resurgence}.

In addition to his direct interventions into mathematical and theoretical physics, Weyl published (and proposed in talks) profound reflections on the epistemology and ontology of the physical world, and  the role of mathematics in it, most notably  \citep{Weyl:PMN,Weyl:Hs91a:31,Weyl:Symmetrie}. Transformations of mathematical structures played a great role in his   reflections; but in stark contrast to Hausdorff  he  proposed to identify, as clearly as possible,  what he considered the 
 automorphisms (global and gauge)  of ``Nature'' herself to which the transformation group of the descriptive symbol system ought to adapt as smoothly as possible.   Weyl's objective-transcendental constructive mode has recently been taken up in the philosophy of physics by \citep{Catren:Klein-Weyl}.  We have seen that Hausdorff used the method of transformation and the related structure groups   with  the opposite  goal of undermining  a belief (at least a naive one)  in being able to discern such structures in the world, i.e., in a deconstructivist mode {\em ante letteram}.

For Weyl,  philosophical reflections seemed  important also for securing  a cultural basis for mathematics, in particular the parts which were not amenable to what he accepted as constructive  (i.e., essentially by denumerable procedures). From the mid-1920s onward he realized, at first  hesitatingly, that the principle of the excluded third and axiomatically postulated transfinite  mathematical objects of higher cardinality  may be of importance and {\em acceptable} because of their  role in making the difficult structures of modern physics intelligible, at least in an indirect symbolic way. 
\begin{quote}
From the formalist standpoint, the transfinite component of the axioms takes the place of complete induction and imprints its stamp upon mathematics. The latter does not consist here of evident truths but is {\em bold theoretical construction},  and as such the very opposite of analytical self-evidence. \ldots

In axiomatic formalism, finally, consciousness makes the attempt to `jump over its own shadow'. to leave behind the stuff of the given, to represent the {\em transcendent }  -- but, how could it be otherwise?, only through the {\em symbol}. \citep[64ff.]{Weyl:PMNEnglish}
\end{quote}
 Hausdorff found joy in searching for logically consistent insight into transfinite constructions in the wide sense; he  considered it  as a goal of its own which  carried an intrinsic value.  Weyl, in contrast, considered such symbolical  thought figures (dealing with a stronger transfinite than denumerable constructivism would accept) as  meaningful only if it could be  related  to  natural sciences  directly or indirectly \citep[61]{Weyl:PMNEnglish}. The difference could not be larger. But does one of these opinion  devaluate the other as a legitimate  position of a 20th century ``modern'' mathematician? We better  consider both as  understandable reactions of creative mathematicians to the challenge of the cultural and social modernization they lived in and contributed to. 

\section{Outlook:  Modernity,  emancipation or  crisis in permance? \label{section modernism}}

\subsection{Hausdorff:  liberation, rationalism and the ``end of metaphysics''  \label{subsection Hausdorff end of metaphysics}}

We are well informed about Hausdorff's perception and evaluation of the cultural development in late 19th century Germany through the publication of his {\em alter ego}  Mongr\'e, in particular his time-critical essays in the {\em Neue Deutsche Rundschau}, a leading journal of liberal intellectuals in Germany  \citep[Chap. 6]{Purkert/Brieskorn:Hausdorff}, \citep[vol. {\em 7}]{Hausdorff:Werke}. As mentioned above, he came from a conservative Jewish  parents home, the religious traditions and creed of which he did not share. He grew  up  in a German environment in rapid modernizing  change, which allowed for a slow and selective emancipation of Jewish people  on the one hand, but on the other hand was also  hatching  a rising anti-Semitism in daily life. On this background Hausdorff  developed a sharp-minded, critical, highly individualistic view of life and culture which at the turn to the 20th century was characterized in Germany by a  streaky mix of turbid tradition and cheered up modernism. Later in his life he characterized his own cultural and philosophical trajectory as having developed 
``from Wagner to Schopenhauer, from there back to Kant and forward to Nietzsche'' \citep[vol. {\em 9}, 503]{Hausdorff:Werke}. With  ``Nietzsche'' Hausdorff  at this place referred to the {\em young} (pre-crisis)  writer  whom he emphatically talked about, at a different  place, as the 
\begin{quote}
\ldots affectionate, tempered, appreciative, freethinking Niezsche and the cool, dogma free, system-less sceptic Nietzsche and the (\ldots) world blessing, all positive ecstatic Zarathustra \citep[181]{Purkert/Brieskorn:Hausdorff}.\footnote{``\ldots von dem g\"utigen,
ma\ss{}vollen, verstehenden Freigeist Nietzsche und von dem k\"uhlen, dogmenfreien, systemlosen Skeptiker Nietzsche und von dem Triumphator des Ja-und Amenliedes, dem weltsegnenden, allbejahenden Ekstatiker Zarathustra'' \citep[181]{Purkert/Brieskorn:Hausdorff}.}
\end{quote}
This picture of Nietzsche stands in stark contrast to  the  later  ``fanatic'' Nietzsche who, in addition,  was  contorted to the worse by his the editors under the leadership of his sister. The late, fanatic Niezsche preached a morality which, according to an  observation made by Hausdorff as early as  1902, contained the potential for  ``turning into a world-historic scandal   which might dwarf the inquisition and the witch trials, such that they would  appear as harmless aberrances''   \citep[180]{Purkert/Brieskorn:Hausdorff}.\footnote{``In Nietzsche gl\"uht ein Fanatiker. Seine Moral der Z\"uchtung, auf unserem heutigen Fundamente biologischen und physiologischen Wissens errichtet: das k\"onnte ein weltgeschichtlicher Skandal weden, gegen den Inquisition und Hexenproze\ss{} zu harmlosen Verirrungen verblassen''  \citep[180]{Purkert/Brieskorn:Hausdorff}.}

In short, the young Hausdorff/Mongr\'e developed into an enlightened  Nietzschean dissident in late 19th century  Germany.  He considered the cultural modernisation as an emancipatory chance, with  the intellectual and social liberation of the individual as the cultural task of the time. Some of his writings as P. Mongr\'e had  the flavour of an,  in  my view (ES),   drastically exaggerated  emphasis on the role of individual perception of the world, and the fiction of the happiness of the  ``higher'' persona  standing  above the happiness of the many and  in contrast to any other kind of social bonds  \citep{Hausdorff:Massenglueck}, \citep[sec. 6.1.1]{Purkert/Brieskorn:Hausdorff}.  To him the motif of individual freedom seemed to fit  well with  Cantor's battle cry for set theory: ``the essence of mathematics is freedom''.    In contrast to Cantor himself,  Hausdorff took set theory  as a chance for  dissolving  thinking  from  metaphysical bonds not only inside mathematics, but in general, with mathematics as a trailblazer. In his view, mathematical thought ought to be tied  back to the social and  outer material world only in the indirect and sceptical form of his ``considerate empiricism'' (see above).

In his book ``The Chaos in Cosmic Selection. An Epistemological Essay'' Hausdorff/Mongr\'e hoped
 to be able to do away  with metaphysics once and for all.  The book ended with the often cited, (all too) proud claim:
 \begin{quote}
 Therewith the bridges have been torn down, which, in the imagination of all metaphysicians, connect  the chaos [the transcendent world, ES] and the cosmos [the ordered sensible and intelligible world, ES] in both directions, and the {\em end of metaphysics} has been declared, the explicit one no less than the masked one, both of which the science of the coming century is obliged to scrap from its architecture \citep[209; {\em 7}, 803, emph. in the original]{Hausdorff:Chaos}.\footnote{``Damit sind die Br\"ucken abgebrochen, die in der Phantasie aller Metaphysiker vom Chaos zum Kosmos her\"uber und hin\"uber f\"uhren, und ist das {\em Ende der Metaphysik } erkl\"art, -- der eingest\"andlichen nicht minder als jener verlarvten, die aus ihrem Gef\"uge auszuscheiden der Naturwissenschaft des n\"achsten Jahrhunderts nicht erspart bleibt'' \citep[209; {\em 7}, 803]{Hausdorff:Chaos}.}
\end{quote}

He broadened the argument in a more popular and widespread article 
  ``The  unclean century'' in the {\em Neue Deutsche Rundschau} \citep{Hausdorff:Jahrhundert}. It  contains a    beautiful, in large parts satirical, general settlement with the cultural inconsistencies of the semi-modern culture in Wilhelmian Germany. Hausdorff/Mongr\'e attacked, among others, the militaristic habitus among the German self-defined elite, still fond of the duel as form of honour-saving conflict resolution, certain aspects of  Neo-Kantianism in German humanistic education as cultural hypocrisy,  the rising  neo-religiousness of diverse flavours as obscurantism, and the unacknowledged metaphysical elements in natural science
 \citep[sec. 6.1.2]{Purkert/Brieskorn:Hausdorff}. In his polemics all this appeared as a hangover of  earlier  times  and had to be done away with by 
\begin{quote}
\ldots an act of cleanliness with which any retiring century should recommend itself to its successor \citep[352]{Hausdorff:Jahrhundert}.\footnote{Man ``\ldots vollzieht einen Act der Reinlichkeit, mit dem jedes scheidende Jahrhundert sich seinem Nachfolger empfehlen sollte'' \citep[Werke {\em 7}, 352]{Hausdorff:Jahrhundert}.}
\end{quote}
For the young Hausdorff (Mongr\'e)    some sort of purified modernity appeared  as a  desirable future state of the human  world. Needless to say that this optimistic perspective was broken by the two Great Wars, the  deep world crisis  of the early 20th century between them,  and the rise of Nazism to power in Germany, with all the humiliations and cruelties against  the Jewish population, which he himself had to go through. One of his last letters written in January 1941, about a year before his enforced suicide, ended with  realistic resignation:
\begin{quote}
Nietzsche always feared that Europe might perish because of a hysteria of pity: one cannot claim that this diagnosis was particularly realistic \citep[vol. {\em 9}, 357]{Hausdorff:Werke}.\footnote{``Nietzsche hat immer bef\"urchtet, dass Europa an einer Hysterie des Mitleids zugrunde gehen w\"urde: man kann nicht behaupten, dass diese Diagnose sehr zutreffend war'' (Hausdorff an J. K\"afer, 2. Jan. 1941,  \citep[vol. {\em 9}, 357]{Hausdorff:Werke})}
\end{quote}

\subsection{Weyl: awareness of crisis  and the search for  metaphysical horizons \label{subsection Weyl metaphysics}}
 Weyl was among those who, while still at school, was strongly affected by Kant's critical philosophy.  For him this did not at all lead to a complacent and indolent attitude, so pungently attacked by Mongr\'e in his essay about the ``unclean century''.  In retrospect he  characterized the effect of  Kant's teaching of the ``ideality of space and time'' quite differently:
 \begin{quote}
 \ldots by one jerk I was awoken from the `dogmatic slumber'; the world was most radically put into question for the mind of the adolescent \citep[{\em 4}, 632]{Weyl:EuB}.\footnote{``\ldots mit einem Ruck war ich aus dem dogmatischen Schlummer' erwacht, war dem Geist des Knaben auf radikale Weise die Welt in Frage gestellt'' \citep[{\em 4}, 632]{Weyl:EuB}.}
 \end{quote}
 Thus, for Weyl,  the reading of Kant had  an effect usually ascribed to ``modernity'' or ``modernism'': a radical detachment of assuming simple bonds to reality. This detachment  was even enhanced, when he entered G\"ottingen university and learned of Hilbert's studies of the foundations of geometry. The  ``multitude of different unfamiliar geometries'' studied in the axiomatic approach  destroyed his simplified picture of an  ``edifice'' of Kantian philosophy, which he had erected in his mind (ibid. 633). This retrospective description indicate that Weyl, in contrast to Hausdorff/Mongr\'e, sensed the  confrontation with a  ``modern'' view of the world, and the adoption of it for himself, as a deeply irritating experience. In much of his later writings we find an embarrassment about the basic  detachment of mathematical knowledge from the link to the external world. Weyl would sometimes speak of a  ``transcendent'' reality, apparently alluding also to the religious connotation of the word besides a  vague reference to an outer nature beyond the one ``given'' to the senses and to phenomenal insight. 

Many authors have argued that   the experience of the breakdown of civil norms during the Great War and the following deep social crisis in Germany aggravated Weyl's, and others, sensitivity with regard to the stability also  of scientific and even mathematical knowledge \citep{Mehrtens:Moderne,Skuli:DMV,Schappacher:Politisches}. The latter had been untightened already in the later 19th century by the loss of credibility of traditional metaphysics and an imputed direct reference to an external reality.  This seems to have strongly influenced Weyl's sensitivity for crisis in the  debate on the foundations of analysis and set theory. 

It seems  that  Weyl experienced the rise of  modernity, i.e.,  of modern society in the sense of late 19th- and 20th-century high capitalism and its scientific culture, as a challenge and a crisis set in permanence. We have seen that 
in the second half of the 1920s he was willing to accept that Hilbert's proof theoretic (``finitist'') program might even be successful with regard to a formal legitimation of the use of (strong) transfinite methods in mathematics. But in his opinion this would not solve the problem of  meaning  of such parts of  mathematics, which were based on transfinite axiomatic methods. In his view (as we know. not in Hilbert's view) G\"odel's incompleteness theorem for a sufficiently strong formal system embracing arithmetic and a formalized logic as strong as the one of Russell's {\em Principia Mathematica} dealt a ``terrific blow'' to Hilbert's  program \citep[{\em 4}, 279]{Weyl:Russell}. This was written  after another, even more devastating war than the one after which he had declared the new ``crisis'' in the  foundations of mathematics. Weyl gave  a short survey of the development of the research in the foundations of mathematics during the last few decades; then he repeated his diagnosis   of the situation, given roughly 30 years ago:
\begin{quote}
From this history one thing should be clear: we are less certain than ever about the ultimate foundations of (logic and) mathematics. Like everybody and everything in the world today, we have our `crisis'. We have had it for nearly fifty years. \citep[{\em 4}, 279]{Weyl:Russell}
\end{quote}
As we also know, this did not hinder him to participate in the enterprise of modern mathematics and physics, but it shaped his selection of research topics and methods. He continued:
\begin{quote}
Outwardly it does not seem to hamper our daily work, and yet I for one confess that it has had a considerable practical influence on my mathematical life: it directed my interests to fields I considered relatively 'safe', and has been a constant drain on the enthusiasm and determination with which I pursued my research work. (ibid.)
\end{quote}
 Remember that `safeness' in the sense of cultural meaning of mathematics included its  link to the clarification of knowledge in the natural science, in particular physics. In addition,  this remark may also be read  as a partial explanation for Weyl's never-ending  efforts to find support in philosophical reflection of his work, an effort which did not stop short of explicit metaphysical considerations.
 This stands in the sharpest possible contrast to    Hausdorff  whose verdict of (classical) metaphysics  we have seen above.

\subsection{Final remarks:  modern -- countermodern -- trans-modern ?  \label{section final}}
Neither of our two protagonists maintained a Platonist view of the objects of mathematical knowledge. Hausdorff rejected any claim of ideal order beyond the  insights gained by  logically precisely framed symbolic production in the realm of transfinite sets opened by Dedekind and Cantor.  He was convinced that such an argumentation can be   expanded without running into contradictions, as long as  a carefully restrained use of the comprehension principle was made. Although Hausdorff did not, to my knowledge, publish a short, conclusive  verbal description of his view of mathematics, he may well be called a {\em symbolic formalist}.
That is, he emphasized that mathematics deals with ``objects of thought, symbols of undetermined meaning''  which underlie no other  constraint than that of  logical consistency.\footnote{The closest approximation to such a short characterisation can be found in section 1, ``Der Formalism'', of  an unpublished fragment \citep{Hausdorff:RuZ}  written about 1904, in particular folio 4ff, vol. {\em 6}, p. 474ff. }   
In this respect he had a  completely different conception of set theory than Cantor (who believed in an ontological meaning of transfinite sets),  one which may rightly  be called ``modernist''. Weyl, as we have seen, had a rather different understanding of what mathematics is about, or  at least ought to be about. His perspective of a constructivist or, over a time period, even  intuitionist understanding (sui generis)  of mathematical objects did not allow him to join the radical modernist  attitudes of Hilbert, Hausdorff and,   later, the young mathematicians of the Bourbaki generation. In its decidedly constructivist perspective it was, however, not at all ``countermodern''. It even had  strong resemblances to certain features of modernist architecture (Bauhaus) or art (cubism). Also Weyl's most important  philosophical inspirations received from   Husserl's phenomenology and  Fichte's ``constructivism'' (as he himself described it in \citep[{\em 4}, 641]{Weyl:EuB} and the latter's contemporary  liberal interpreter Fritz Medicus cannot be  qualified as a ``countermodern'' (in contrast to the conservative nationalist interpreters of the South-German Fichteans), or even as an  ``antimodern''  influence on Weyl.  

Finally, if one takes the corpora of the mathematical research work into account, surely the most important sources for the description of a mathematicians, we see here two towering figures of mathematics in the ``modern'' period  of the late 19th and the 20th century. It would be beside the track  to describe one of them, Weyl, as a ``countermodern'' mathematician and only the other one, Hausdorff, as modern. But, of course, the  qualification of Weyl's and Brouwer's position in the foundations of mathematics as representatives of  ``countermodernism'' (``Gegenmoderne'') in Herbert Mehrtens influential book  \citep[289ff., 301]{Mehrtens:Moderne} is not  without any factual base. It  resides on  real  differences between the two authors, which may be described in simple terms as follows: 
While Hausdorff, at least as a young man, welcomed the rising modernity/modernism in science and culture enthusiastically as a liberating movement, Weyl was irritated and suffered from the loss of security and kept distance to modernist positions in the reflective disourse on mathematics (as Mehrtens calls it). This made him  a modern scientist (not a modernist), who was critical of many aspects of modernity, not only with regard to epistemic questions but also with regard to the social destructions which were part of  the rise of modernity. 

 After the second Great War of the 20th century he was shocked by the destructive potential which had developed on the basis of scientific achievements. 
In a manuscript written close to the end of his life and  published  posthumously by T. Tonietti, he deplored the state of things,  dramatically clothed in his grave,  humanistic style. He warned that modern science may be characterized by a kind of {\em hubris} (violent arrogance)  and went on 
 \begin{quote}
For who can
close his eyes against the menace of our own self-destruction by science? The alarming fact is that the rapid progress of scientific knowledge is not paralleled by a
corresponding growth of man's moral strength and responsibility, which have hardly changed in historical
times  \citep[12]{Weyl:1985}.
\end{quote}
Weyl could not even follow Hardy's  move   after the first  War for  exculpating pure mathematics on the basis of  its ``uselessness'' in practical matters, which, according to Hardy, would protect it against a  participation in ``exploitation of our fellow-men''  and the destruction up to their ``extermination'' (Weyls words). Weyl did not believe in such a  escape route and emphasized: 
\begin{quote}
However the
power of science rests on the combination of experiment, i.e., observation under freely chosen conditions,
with symbolic construction, and the latter is its mathematical aspect. Thus if science is found guilty, mathematics cannot evade the verdict.
(ibid.)
\end{quote}
We remember that a  similar, although less dramatically stated
fright was expressed in  Hausdorff's  downplayed  remark of 1941 on Nietzsche's not ``particularly realistic'' warning that modern history might suffer from too much pity for the fellow-men or nature. 

Weyl was a sceptic  modern actor all over his life. As we  know too well, the dangers  of extermination of mankind by war and/or destruction of our natural habitat are now  even more severe than in the 1950s. But science is not only an accomplice of the destructive sides of modernity; it  also plays the role of collecting the warning signs and is necessary for exploring exit strategies from the ongoing destruction.  From our vantage point of the early 21st century,  Weyl  may appear  as a  modern scientist who  tried to dive through the wave of modernism towards some not yet clearly visible type of  trans-modern culture. The latter would mean to stay true to the enlightened elements of modernity, but   to get rid of the destructive forces against  nature and  our fellow-men. 

\subsubsection*{Acknowledgements:}I thank Walter Purkert, Norbert Schappacher, Jos\'e Ferreir\'os, and Lizhen Ji for  their detailed valuable comments to the   preliminary version of this chapter.

\end{document}